\documentclass{amsart}

\normalsize

\newtheorem{theorem}{Theorem}[section]
\newtheorem{lemma}[theorem]{Lemma}
\newtheorem{proposition}[theorem]{Proposition}

\theoremstyle{definition}

\theoremstyle{remark}

\numberwithin{equation}{section}

\newcommand{\N}{\mathbf{N}}
\newcommand{\Z}{\mathbf{Z}}

\newcommand{\ko}{\: , \;}

\renewcommand{\tilde}[1]{\widetilde{#1}}

\newcommand{\ra}{\rightarrow}

\newcommand{\da}{\downarrow}
\newcommand{\ua}{\uparrow}

\newcommand{\arr}[1]{\stackrel{#1}{\rightarrow}}

\newcommand{\iso}{\stackrel{_\sim}{\rightarrow}}

\newcommand{\opname}[1]{\operatorname{\mathsf{#1}}}

\newcommand{\Mod}{\opname{Mod}\nolimits}

\newcommand{\id}{\mathbf{1}}
\newcommand{\R}{\mathbf{R}}
\renewcommand{\L}{\mathbf{L}}

\newcommand{\Bimod}[2]{\opname{Bimod}\nolimits}
\newcommand{\Shmod}[2]{\opname{Shmod}\nolimits}

\newcommand{\op}{^{op}}
\newcommand{\ten}{\otimes}
\newcommand{\tp}[1]{^{\ten #1}}

\newcommand{\im}{\opname{im}\nolimits}

\newcommand{\Hs}{{H}^*}

\newcommand{\cb}{{\mathcal B}}
\newcommand{\cc}{{\mathcal C}}
\newcommand{\cd}{{\mathcal D}}

\newcommand{\ch}{{\mathcal H}}

\newcommand{\eps}{\varepsilon}

\renewcommand{\phi}{\varphi}

\newcommand{\Hom}{\opname{Hom}}
\newcommand{\Homdot}{\opname{Hom}^\bullet}

\newcommand{\Hominf}{\raisebox{0ex}[2ex][0ex]{$\overset{\,\infty}{
                              \raisebox{0ex}[1ex][0ex]{$\mathsf{Hom}$}
                                                                 }$}}
\newcommand{\Hominfdot}{\Hominf^\bullet}

\newcommand{\Tot}{\opname{Tot}}

\begin{document}

\title{Bimodule complexes via strong homotopy actions}
%\dedicatory{To Professor Klaus W. Roggenkamp on the occasion of his
%sixtieth birthday}
\author{Bernhard Keller}
\address{UFR de Math\'ematiques\\
   UMR 7586 du CNRS \\
   Case 7012\\
   Universit\'e Paris 7\\
   2, place Jussieu\\
   75251 Paris Cedex 05\\
   France }

\email{
\begin{minipage}[t]{5cm}
keller@math.jussieu.fr \\
www.math.jussieu.fr/ $\tilde{ }$ keller
\end{minipage}
}

\date{September 20, 1999}

\begin{abstract}
We present a new and explicit method for lifting
a tilting complex to a bimodule complex. 
The key ingredient of our method is the notion of a strong
homotopy action in the sense of Stasheff.
\end{abstract}

\maketitle

%\tableofcontents

\section{Introduction}

Let $A$ and $B$ be (associative, unital) algebras over
a commutative ring $k$. 
Denote by $\Mod A$ the category of
(right) $A$-modules. Suppose that $P$ is a $B$-module
endowed with an algebra morphism
\[
A \ra \Hom_{\Mod B}(P,P).
\]
Then $P$ becomes an $A$-$B$-bimodule and we have the
tensor functor
\[
?\ten_A P : \Mod A \ra \Mod B\ko
\]
which takes the free $A$-module $A_A$ to $P$.
This is the basic fact which allows us to construct
Morita equivalences $\Mod A \iso \Mod B$.

Now let $\cd B=\cd \Mod B$ be the (unbounded) derived
category of the category of $B$-modules. Suppose that we
have a complex $T\in \cd B$ and an algebra morphism
\[
A \ra \Hom_{\cd B}(T, T).
\]
It is not clear whether this map comes from an action
of $A$ on the \emph{components} of $T$, even
after replacing $T$ by an isomorphic object of $\cd B$.
Therefore, the (derived) tensor product by $T$
not well-defined and the analogy with the case
of module categories seems to break down.
Nevertheless, for complexes $T$ satisfying the `Toda condition'
\[
\Hom_{\cd B}(T, T[-n])=0\ko \mbox{ for all }n>0\ko
\]
J.~Rickard succeeded \cite{Rickard89a} in constructing
a triangle functor $F:\cd^- A \ra \cd B$ taking $A_A$ to $T$,
where $\cd ^{-}A\subset \cd A$ denotes the subcategory of 
right bounded complexes.
This construction was at the basis of the proof of his
`Morita theorem for derived categories'. Later he showed
\cite{Rickard91} that if $F$ restricts to an equivalence 
between the bounded derived categories 
(and suitable flatness hypotheses
hold), then after replacing $T$ by an isomorphic object,
it \emph{is possible} to lift the $A$-action to the components of $T$.
In his proof, he used the functor $F$ constructed in \cite{Rickard89a}.

In \cite{Keller93} and \cite{Keller94}, we gave an a priori
construction of a lift (up to isomorphism) of $T$ to a
complex of bimodules $X$ (under suitable flatness hypotheses). 
This made it possible to \emph{define}
$F=\L(?\ten_A X)$ and to give a new proof
\cite[Ch. 8]{KoenigZimmermann98} of Rickard's
Morita theorem. 

Thus, up to now, there have been two constructions of a bimodule
complex $X$ from a complex $T$ as above. Neither of them is
very explicit: the first one \cite{Rickard91} uses the
functor $F$; the second one \cite{Keller93} uses resolutions
over differential graded algebras.

In this paper, we present a new construction, which is 
surprisingly explicit. In fact, if we assume that $T$ is a right
bounded complex of projective $B$-modules, then essentially
the only data we need are homotopies $H(f)$ such that
\[
f= d\circ H(f) + H(f)\circ d
\]
for each morphism of complexes $f:T \ra T[-n]$, $n>0$.
We also prove a unicity result which improves on \cite{Rickard91}
and \cite{Keller93}. 

The essential new ingredient of our method is the notion
of a \emph{strong homotopy action} (=$A_\infty$-action) due to
Stasheff \cite{Stasheff63a}, \cite{Stasheff63b} and
recently popularized again by Kontsevich \cite{Kontsevich94},
\cite{Kontsevich98}. The present article is self-contained
but the interested reader may find more information on 
strong homotopy methods in \cite{Keller99a}. 
We will show that if $A$ is projective over $k$ and $T$
is right bounded with projective components and satisfies
the Toda condition, then the `action up to homotopy'
of $A$ on $T$ may be enriched to a strong homotopy
action. It is remarkable that this can be done without changing
the underlying complex of $T$. In a second step, we
show that each strong homotopy action on a complex
$K$ yields a strict action on a larger (but quasi-isomorphic)
complex $K'$. In fact, $K'$ may be viewed as a 
`perturbed Hochschild resolution' of the complex $K$.

\bigskip
\emph{Acknowledgment.} The author wishes to thank 
Prof.~K.~W.~Roggenkamp for
his stimulating interest and encouragment.

The work on this article was begun during a stay
of the author at the Sonderforschungsbereich `Diskrete Strukturen in der 
Mathematik' at the University of Bielefeld.
The author thanks the Bielefeld representation theory 
group, and in particular C.~M.~Ringel
and P.~Draexler, for their kind hospitality.

\tableofcontents

\section{The main theorem}\label{main}

Let $k$ be a commutative ring and $B$ an (associative,
unital) $k$-algebra. Let $T$ be a complex
of (right) $B$-modules. Let $A$ be another $k$-algebra. 
A \emph{strict action} of $A$ on $T$ is an homomorphism
(preserving the unit)
\[
A \ra \Hom_{\cc B}(T,T)\ko
\]
where $\cc B$ is the category of complexes of
$B$-modules. Equivalently, a strict action is the datum of a complex
of $A$-$B$-bimodules whose restriction to $B$ equals $T$.
An \emph{homotopy action} of $A$ on $T$ is an homomorphism
\[
\alpha : A \ra \Hom_{\ch B}(T,T)\ko 
\]
where $\ch B$ is the homotopy category of right $B$-modules.
The pair $(T,\alpha)$ is then an \emph{homotopy module}.
If $T$ and $T'$ are endowed with homotopy actions $\alpha $ and $\alpha '$,
a morphism of complexes $f:T\ra T'$ is \emph{compatible} with these
actions if $f\circ \alpha(a)$ is homotopic to $\alpha'(a)\circ f$
for all $a\in A$. Then $f$ is also called a \emph{morphism
of homotopy modules}.

Clearly, each strict action yields an homotopy action.
The converse is false, in general. However, we have the

\begin{theorem}\label{maintheorem} Let $k$ be a commutative ring and
$A$, $B$ two (associative, unital) $k$-algebras. Let $T$ be a complex
of right $B$-modules endowed with an homotopy action by $A$. Suppose
that $A$ is projective as a $k$-module, and that $T$ is right bounded
with projective components and satisfies
\begin{equation}\label{Toda}
\Hom _{\ch B}(T,T[-n])=0 \text{ for all } n>0.
\end{equation}
\begin{itemize}
\item [a)] There is a right bounded complex of projective $B$-modules
$X$ endowed with a strict action by $A$ and a quasi-isomor\-phism
$\phi : T\ra X$ of complexes of $B$-modules compatible with the
homotopy actions by $A$.
\item[b)] If $\phi:T\ra X$ and $\phi' : T \ra X'$ are two
quasi-isomorphisms as in a), then there is a unique isomorphism
$\psi:X \ra X'$ in the derived category of $A$-$B$-bimodules such that
we have $\psi\circ \phi =\phi'$ in the homotopy category of
$B$-modules.
\end{itemize}
\end{theorem} 
For the case where $T$ is a tilting complex, the theorem was proved by
J.~Rickard in \cite{Rickard91}.  The general case is proved in
\cite{Keller93}. We give a new proof which, compared to these previous
approaches, is much more explicit. We illustrate this by two special
cases.

\section{First special case} \label{firstspecial}

Suppose that the assumptions of theorem \ref{maintheorem} hold.
Let us assume that  $T$ has non vanishing components at most in degrees 
$0$ and $1$:
\[
T=(\dots \ra 0 \ra T_{1}\ra  T_{0}\ra 0\ra \dots ).
\]
We will construct a quasi-isomorphism $\phi :T \ra X$ as in theorem
\ref{maintheorem} a) \emph{without using property \eqref{Toda}}.

Since $A$ is projective over $k$, we can find a $k$-linear map
$\tilde{\alpha }: A \ra \Hom _{\cc }\, (T,T)$ lifting the
given homotopy action $\alpha :A \ra \Hom _{\ch B}\,(T,T)$.
Now we define $m_{2}:A\ten T \ra T$ by
\[
m_{2}(a,x)=(\tilde{\alpha }(a))(x)\ko a\in A, x\in T.
\]
Now we have to take into account the non-associativity
of $m_2$: Consider the square
\[
\begin{array}{rcl}
A\ten A\ten T & \arr{\id_A \ten m_{2}} & A\ten T \\
m_A\ten \id_T\da  &                     & \da m_{2} \\
A\ten T &  \arr{m_{2}}      & T\ko
\end{array}
\]
where $m_A$ denotes the multiplication of $A$. 
The square becomes commutative in the homotopy category. Hence there is
a morphism of graded $A$-modules 
\[
m_{3}:A\ten A\ten T \ra T
\]
homogeneous of degree $-1$ such that we have
\[
m_{2}(ab,x) - m_{2}(a,m_{2}(b,x)) = -m_{3}(a,b,d(x))-d(m_{3}(a,b,x))
\]
for all $a,b\in A$ and $x\in T$.
We construct a complex $\tilde{X}$ as follows: 
The underlying graded module of $\tilde{X}$ is 
\[
A\ten A\ten A\ten T[2]\;\; \oplus \;\;A\ten A\ten T[1] \;\;
\oplus \;\; A\ten T.
\]
The differential is given by
\begin{align*}
d(a,b,c,x) & = -ab\ten c\ten x + a\ten bc \ten x \\
           & \quad + a\ten m_3(b,c,x) - a\ten b\ten m_2(c,x)
               + a\ten b\ten c\ten d(x)\ko \\
d(a,b,x) & =  -ab\ten x + a\ten m_{2}(b,x)- a\ten b \ten d(x)\ko \\
d(a,x)   & =  a\ten d(x)\ko
\end{align*}
where $a,b,c\in A$, $x\in T$.
We define the complex $X$ to be the truncation
\[
X=\tau_{\leq 1} \tilde{X}= (\tilde{X}_1/\im {d_2} \ra \tilde{X}_0).
\]
We define the morphism of complexes $f: T \ra \tilde{X}$ 
by $f(x)=1\ten x$. This morphism is compatible with the
homotopy actions by $A$. Indeed, if we define the graded
morphism $f_2: A \ten T \ra \tilde{X}$ of degree $-1$ by
\[
f_2(a\ten x) = 1\ten a\ten x \in A\ten A\ten T[1]\ko a\in A, x\in T\ko
\]
then we have
\[
d(f_{2}(a,x))+ f_{2}(a,d(x)) =
f(m_{2}(a,x))-m_{2}(a,f(x)) \ko a\in A\ko x\in T.
\]
By composition, $f$ yields a morphism $T \ra X$.
Its homotopy class is the required quasi-isomorphism 
$\phi$. Since $T \ra \tilde{X}$ is compatible with the
homotopy action by $A$ and $\tilde{X} \ra X$ is a
morphism of complexes of bimodules, the morphism $\phi$ is
compatible with the homotopy action by $A$.
Note that we have not used the vanishing property \eqref{Toda} of
the complex $T$.

\section{Second special case}\label{secondspecial}

Suppose that the assumptions of theorem \ref{maintheorem} hold.
Let us assume that  $T$ has
non vanishing components at most in degrees $0$, $1$ and $2$:
\[
T=(\dots \ra 0 \ra T_{2}\ra T_{1}\ra  T_{0}\ra 0\ra \dots ).
\]
We will construct a complex of bimodules $X$ and a quasi-isomorphism
$\phi :X \ra T$ as in theorem \ref{maintheorem} a).
For this, we construct a morphisms $m_2$, $m_3$
as in section \ref{firstspecial}. Now consider
the graded morphism $c:A^{\ten 3}\ten T \ra T$ defined by
\[
c(a,b,c,x)=
m_{3}(ab,c,x) -m_{3}(a,bc,x)+m_{3}(a,b,m_{2}(c,x))-m_{2}(a,m_{3}(b,c,x)).
\]
A computation shows that $c$ defines a morphism of complexes
$A^{\ten3}\ten T \ra T[-1]$. By the \eqref{Toda}, there
exists a graded morphism $m_{4}:A^{\ten 3}\ten T \ra T$ homogeneous
of degree $-2$ such that we have
\[
c(a,b,c,x)= m_{4}(a,b,c,d_{T}(x)) - d_{T}(m_{4}(a,b,c,x))\ko
\]
$a,b,c\in A$, $x\in T$. We define $\tilde{X}$ to be
the complex whose underlying graded module is
\[
\bigoplus_{i=1}^4 A\tp{i} \ten T[i-1]
\]
and whose differential is given by
\begin{align*}
d(a_0,\ldots, a_3, x) &= - (a_0 a_1, a_2, a_3, x)
                         + (a_0, a_1 a_2, a_3, x) - (a_0, a_1, a_2 a_3, x)\\
                      &\quad + (a_0, m_4(a_0, a_1, a_2, a_3, x))
                         - (a_0,a_1, m_3 (a_2, a_3, x)) \\
                      &\quad + (a_0,a_1,a_2, m_2(a_3,x))
                         - (a_0, a_1, a_2, a_3, d(x)) \\
d(a_0, a_1, a_2,x) &= -(a_0 a_1, a_2, x) + (a_0, a_1 a_2, x) \\
                   &\quad + (a_0, m_3(a_1,a_2,x)) - (a_0,a_1, m_2(a_2,x))
                          + (a_1, a_1, a_2, d(x)) \ko \\
d(a_0, a_1, x) &= -(a_0 a_1,x) + (a_0, m_2(a_1,x)) \ko \\
d(a_0,x)   &= (a_0, d(x)). 
\end{align*}
We define $f: T \ra \tilde{X}$ by $x \mapsto 1\ten x$ and
we define $X$ to be the truncation
\[
X=\tau _{\leq 2}\tilde{X}.
\]
The homotopy class of the composition $T \ra \tilde{X} \ra X$
is the required morphism $\phi$. As in section \ref{firstspecial},
one checks that $\phi$ is compatible with the homotopy actions
by $A$.

\section{Proof of unicity}\label{proofunicity}

We will prove part b) of the main theorem. This 
could be done by strong homotopy methods as well. 
The following argument is shorter but less explicit.

Let us first observe that $\phi $ and $\phi '$ are homotopy
equivalences of complexes of $B$-modules. 
So there is a unique morphism $f:X\ra X'$ of $\ch B$ such
that $f\circ \phi=\phi'$ in the homotopy category of $B$-modules.
Of course, $f$ is a morphism of homotopy modules. 
We have to show that it lifts to a unique morphism $X \ra X'$ in 
the derived category of $A$-$B$-bimodules $\cd (A\op\ten B)$. 
Let us compute $\Hom _{\cd (A\op\ten B)}(X,X')$. 
The complex $X$ is quasi-isomorphic to
the (sum) total complex of its Hochschild resolution:
\[
\ldots \ra A\ten A\tp{p} \ten X \ra \ldots \ra A\ten X \ra 0 \ko
p\geq 0.
\]
This total complex is right bounded and its components
are projective over $A\op\ten B$ since $A$ is projective over
$k$ and the components of $X$ are projective over $B$.
So we can compute $\R \Homdot_{A-B}(X,X')$ by applying
$\Homdot_{A-B}(?,X')$ to the total complex of the 
Hochschild resolution.  Using the isomorphism
\[
\Homdot_{A-B}(A\ten A\tp{p}\ten X,X')=\Homdot_{B}(A\tp{p}\ten X,X')
\]
we find that $\R \Homdot_{A-B}(X,X')$ isomorphic to the 
product total complex of the following double complex $D$
\[
\Homdot_B\,(X,X') \ra \Homdot_B\,(A\ten X, X') \ra 
\Homdot_B\,(A\ten A\ten X, X') \ra \ldots .
\]
We have to to compute $H^{0} \Tot ^{\Pi }D$.
For this, we first truncate the columns of $D$: For a complex
of $k$-modules $K$, let 
\[
\tau _{\geq 0}K=(\dots \ra 0 \ra 0\ra K^{0}/\im d^{-1} 
\ra  K^{1}\ra K^{2}\ra \dots ).
\]
Let $D_{\geq 0}$ be the double complex obtained by applying
$\tau _{\geq 0}$ to each column of $D$ and let $D_{<0}$
be the kernel of $D\ra D_{\geq 0}$. We claim that $D_{<0}$ is
acyclic. Indeed, the homology of the $p$-th column of $D_{<0}$
in degree $-q$ is isomorphic to 
\[
\Hom_{\ch\cb}(A\tp{p}\ten X, X'[-q]).
\]
This vanishes for $-q<0$ by the projectivity of $A$ over $k$ and
assumption \eqref{Toda}. Hence each column of $D_{<0}$ is acyclic.
Moreover, $D_{<0}$ is concentrated in the right half plane.
We claim that the product total complex $\Tot ^{\Pi }D_{<0}$ 
is acyclic. Indeed, this complex is the inverse limit of
the sequence of the complexes $A_{p}=\Tot ^{\Pi }F_{p}D_{<0}$, 
$p\geq 0$, associated with the column filtration $F_{p}D_{<0}$.
The $A_{p}$ are acyclic by induction on $p$. Each map
$A_{p+1}\ra  A_{p}$ is surjective in each component. 
It follows that it is surjective in the boundaries and hence in
the cycles (which equal the boundaries since the $A_{p}$ are
acyclic). Therefore the inverse limit of the $A_{p}$ is acyclic,
i.e. $D_{<0}$ is acyclic. So the morphism $D\ra D_{\geq 0}$
induces a quasi-isomorphism in the product total complexes.
Hence it is enough to compute $H^{0}\Tot ^{\Pi }D_{\geq 0}$.
It is straightforward to check that this group is canonically
isomorphic to the group of morphisms of homotopy modules
$X \ra X'$.

\section{From homotopy actions to strong homotopy actions}
\label{fromhomotopyactions}

\subsection{Lifting homotopy actions}\label{liftinghtpyactions}
Suppose that $k$ is a commutative ring, $A$ and $B$ are
associative unital $k$-algebras, and $L$ is a $\Z$-graded $B$-module. 
A \emph{strong homotopy action}
of $A$ on $L$ is the datum of graded ($B$-linear) morphisms
\[
m_{n}:A\tp{n-1} \ten L \ra L
\]
defined for $n\geq 1$ and homogeneous of degree $2-n$ such that for
each $n\geq 1$ and all $a_i\in A$, $x\in L$, we have 
\begin{equation}\label{ainfmodule}\begin{split} 
0 &= \sum _{i=1}^{n-2} (-1)^{i-1} 
         m_{n-1}(a_{1},\dots ,a_{i}a_{i+1},\dots ,a_{n-1},x)+ \\
  &  \quad \sum _{k=1}^{n} (-1)^{n-k} 
m_{n-k+1}(a_{1},\dots ,a_{n-k}, m_{k}(a_{n-k+1},\dots ,a_{n-1},x)).
\end{split}
\end{equation}
Note that if $L$ has non vanishing components only in degrees
$0$, \dots , $l$, then $m_{n}$ vanishes for $n>l+2$. It is instructive
to consider the cases $n=1,2,3$ of \eqref{ainfmodule}: For
$n=1$, we obtain
\[
0=m_{1}m_{1}
\]
so that $(L,m_{1})$ is a complex. For $n=2$, we have
\[
0=-m_{2}(a,m_{1}(x))+m_{1}(m_{2}(a,x))\ko a\in A, x\in L\ko
\]
so that $x\mapsto m_{2}(a,x)$ is a morphism of complexes for all $a\in A$.
For $n=3$, $a,b\in A$, $x\in L$, we have
\begin{equation} \label{m2m3}
0=m_{3}(a,b,m_{1}(x))+m_{2}(ab,x)-m_{2}(a,m_{2}(b,x))+m_{1}(m_{3}(a,b,x))\ko
\end{equation}
which expresses the fact that $m_{2}$ is an associative operation up to 
an homotopy given by $m_{3}$.

\begin{theorem}\label{exstrongaction} Suppose that $L$ is a
graded $B$-module endowed with three graded morphisms
$m_i: A\tp{i-1}\ten L \ra L$, $1\leq i \leq 3$, homogeneous
of degree $2-i$ and satisfying \eqref{ainfmodule} for $n\leq 3$. 
Suppose that we have
\begin{equation} \label{Todastrongaction}
\Hom_{\ch B}\,(A\tp{n}\ten L,L[2-n])=0 \quad \mbox{ for all }n\geq 3.
\end{equation}
Then the triple $m_{1},m_{2}, m_{3}$ may
be completed to a strong homotopy action $m_{n}, n\geq 1$, of
$A$ on $L$.
\end{theorem}

In the next subsection, we will set up the dictionary between
strong homotopy and differential coalgebra. We will then prove the theorem
in \ref{proofexstrongaction} using this dictionary.

\subsection{Differential coalgebra}\label{differentialcoalgebra}
Passing from strong homotopy notions to differential coalgebra
notions is a classical device,
cf. \cite{Stasheff63a}, \cite{Stasheff63b}, \cite{Kadeishvili85}.
In this subsection, we adapt it to our needs.

Let 
\[
C=T(A[1])= k\oplus A[1]\oplus (A[1]\ten A[1]) \oplus \dots 
\]
be the graded tensor algebra over the graded $k$-module $A[1]$.
It becomes a graded coalgebra for the comultiplication defined by 
\[
\begin{split} 
\Delta (a_{1},\dots ,a_{n})= 
1\ten (a_{1},\dots ,a_{n})+ & a_{1}\ten (a_{2},\dots ,a_{n}) + \dots \\
       & + (a_{1},\dots ,a_{n-1})\ten a_{n} +(a_{1},\dots ,a_{n})\ten 1.
\end{split}
\]

The graded coalgebra $T(A[1])$ admits a unique graded endomorphism
$b$ of degree $+1$ which satisfies 
\begin{align*}
b(x_{1},x_{2}) &= x_{1}x_{2} \mbox{ for } x_{1},x_{2}\in A\ko\\
b(x_{1},\ldots, x_{i}) &= 0  \mbox{ for } i\neq 2\ko
\end{align*}
and which is a coderivation:
\[
\Delta \circ b = (b\ten \id_C + \id_C \ten b)\circ \Delta .
\]
Here and elsewhere, we use the \emph{graded tensor product}: 
For graded maps $f$, $g$ and homogeneous elements $x$, $y$, 
we have
\[
(f\ten g)(x\ten y)=(-1)^{|g||x|}f(x)\ten g(y).
\]
where the bars indicate the degree.
Explicitly, the formula for $b$ is
\[
b(x_{1},\dots ,x_{n})=\sum _{i=1}^{n-1} (-1)^{i-1} 
(x_{1},\dots ,x_{i}x_{i+1},\dots ,x_{n})\ko n\geq 2.
\]
We have $b^2=0$ so that $(C,d)$ is a differential graded coalgebra.

Now let $L$ be a $\Z$-graded $B$-module.
Let $X$ be the graded $B$-module $C\ten L$. So we have
\[
X=L \oplus (A[1]\ten L) \oplus \dots \oplus 
(A[1]\tp{i}\ten L) \oplus \cdots \ko i\in \N .
\]
The graded module $X=C\ten L$ becomes a cofree graded comodule 
over $C$ for the comultiplication induced from that of $C$. So
we have
\[
\begin{split}
\delta (a_{1},\dots ,a_{n-1},x)= 1\ten (a_{1},\dots ,a_{n-1},x) 
+ & a_{1}\ten (a_{2},\dots ,a_{n-1},x) + \dots \\
  & + (a_{1},\dots ,a_{n-1})\ten x.
\end{split}
\]
A \emph{coderivation} of $X$ of degree $e$ is a graded endomorphism
$b$ of degree $e$ such that
\[
\delta \circ b = b_C\ten \id_X + \id_C \ten b.
\]
Let $\eps: X\ra L$ be the projection. Then the map $b\mapsto \eps\circ b$
is a bijection from the set of degree $e$ coderivations of $X$ to the set
of graded morphisms of degree $e$ from $X$ to $L$. Let us describe the
inverse map: Let $b: X \ra L$ be a graded morphism of degree $e$.
It is given by its components (homogeneous of degree $e$)
\[
b_i : A[1]\tp{i-1} \ten L \ra L\ko i\geq 1.
\]
The corresponding coderivation is given by
\begin{equation}\label{diffcomodule}\begin{split} 
b(a_1,\ldots, a_{n-1},x) &= \quad \sum _{i=1}^{n-2} (-1)^{e(i-1)} 
               (a_{1},\dots ,a_{i}a_{i+1},\dots ,a_{n-1},x) \\
  &\quad + \sum _{k=1}^{n} (-1)^{e(n-k)} 
       (a_{1},\dots ,a_{n-k}, b_{k}(a_{n-k+1},\dots ,a_{n-1},x))\ko
\end{split}
\end{equation}
where $a_i\in A$, $x\in L$.

Let $U$ be a $C$-comodule. Since $X$ is cofree, the map
$f\mapsto \eps\circ f$ is a bijection from the degree $e$
comodule morphisms $f: U \ra X$ to the degree $e$ graded
$k$-linear morphisms $U \ra L$.

Let $b$ be a degree $1$ coderivation of $C\ten X$. Then we have
\[
(\id_C\ten b)(b_C\ten \id_X)= -(b_C\ten \id_X)(\id_C\ten b)
\]
as morphisms $C\ten X \ra C\ten X$.
Since $b_C^2=0$, we deduce that
\[
b^2 \circ \delta = (\id_C\ten b + b_C \ten \id_X)^2 \circ \delta
= (\id_C \ten b^2)\circ \delta.
\]
It follows that $b^2: X \ra X$ is a morphism of $C$-comodules.
In particular, we have $b^2=0$ iff $\eps\circ b^2=0$. Thanks
to \eqref{diffcomodule}, this last equality translates into
\begin{equation}\label{diffsquarecomodule}\begin{split} 
0 = & \quad \sum _{i=1}^{n-2} (-1)^{(i-1)} 
               b_{n-1}(a_{1},\dots ,a_{i}a_{i+1},\dots ,a_{n-1},x) \\
  & + \sum _{k=1}^{n} (-1)^{(n-k)} 
   b_{n-k+1}(a_{1},\dots ,a_{n-k}, b_{k}(a_{n-k+1},\dots ,a_{n-1},x))\ko
\end{split}
\end{equation}
for all $n\geq 1$.
If we compare  \eqref{ainfmodule} to \eqref{diffsquarecomodule},
we see that the map $(m_n)\mapsto (b_n)$ defined by
\[
b_n(a_1, \ldots, a_{n-1}, x )= m_n(a_1, \ldots, a_{n-1},x)
\]
is a bijection between strong homotopy actions of $A$ on $L$
and degree $1$ comodule differentials on $X$. Note that
$m_i$ is a graded map $A\tp{i-1}\ten L \ra L$ of degree $2-i$
whereas $b_i$ is a graded map $A[1]\tp{i-1}\ten L \ra L$
of degree $1$.

Now assume that $b$ is a degree $1$ coderivation of
$X$. Let us analyse the equation $b^2=0$ in terms
of the components $b_i : A[1]\tp{i-1}\ten L \ra L$.
If we identify them with their 
extensions to coderivations, then the equation $b^2=0$ 
translates into
\begin{equation}\label{component}
0=b_1 b_n + b_2 b_{n-1} + \ldots + b_n b_1 
\end{equation}
for all $n\geq 1$.
Note that for each $p\geq 0$, the right hand side takes 
$A[1]\tp{n+p}\ten L$ to $A[1]\tp{p}\ten L$.
The equation $b^2=0$ holds iff, for all $n\geq 1$,
the right hand side of
\eqref{component} induces the zero map $A[1]\tp{n-1}\ten L\ra L$.
Indeed, in this case we have $\eps\circ b^2=0$.

\subsection{Proof of theorem \ref{exstrongaction}}
\label{proofexstrongaction}
We use the notations of the previous two subsections.
If we identify graded maps $A\tp{i-1}\ten L \ra L$ with their
extensions to coderivations $X \ra X$, then we are given
coderivations $b_1$, $b_2$, $b_3$ such that we have
\begin{align*}
0 =& b_1^2 \ko\\
0 =& b_1 b_2 + b_2 b_1 \ko\\
0 =& b_1 b_3 + b_2 b_2 + b_3 b_1.
\end{align*}
We have to construct $b_i$, $i\geq 4$, such that
\begin{equation} \label{bicondition}
b_1 b_n + b_2 b_{n-1} + \cdots b_{n-1} b_2 + b_n b_1 =0
\end{equation}
for all $n\geq 4$. Suppose that $N\geq 4$ and that
$b_1, \ldots, b_{N-1}$ have been constructed such that 
\eqref{bicondition} holds for all $n\leq N-1$. 
We are looking for $b_N$ such that
\[
0= b_1 b_N + b_N b_1 + (b_2 b_{N-1} + b_3 b_{N-2} + \cdots + b_{N-1} b_2).
\]
Put 
\[
c=b_2 b_{N-1} + b_3 b_{N-2} + \cdots + b_{N-1} b_2.
\]
Let $X_N\subset X$ be the $C$-subcomodule
\[
L\oplus (A[1]\ten L) \oplus \ldots \oplus (A[1]\tp{N-1}\ten L).
\]
Note that $b_i$ takes $X_N$ to $X_{N-i+1}$. In particular,
$c$ takes $X_N$ to $X_1=L$ and vanishes on $X_{N-1}$. So it
induces a graded morphism of degree $2$
\[
X_N/X_{N-1}=A[1]\tp{N-1}\ten L \ra L.
\]
We only have to show that this is a morphism of complexes:
Indeed, by our assumption on $L$, it will then have to be nullhomotopic. 
The extension of an homotopy to a coderivation yields the required 
morphism $b_N$.  
To show that the morphism induced by $c$ commutes with the 
differential $b_1$, we define $B=b_1 + \cdots + b_{N-1}$. 
Since \eqref{bicondition} holds for $n\leq N-1$, we have
\begin{equation}\label{Bsqeqc}
B^2\equiv_N c
\end{equation}
where $\equiv_N$ denotes the equality of the restrictions to $X_N$.
We conclude that $B^2$ vanishes on $X_{N-1}$ and takes $X_N$ to 
$L\subset X_N$. This implies that $B^2 B \equiv_N B^2 b_1$ and 
$B B^2 \equiv_N b_1 B^2$. 
Therefore we have $B^2 b_1 \equiv_N b_1 B^2$ and 
finally $c b_1 \equiv_N b_1 c$.

\subsection{Morphisms of strong homotopy actions}\label{morphismsstronghtpy}

Suppose that $k$ is a commutative ring, $A$ and $B$ are
associative unital $k$-algebras, and $L$ and $M$ are
\emph{strong homotopy modules}, i.e. 
$\Z$-graded $B$-modules endowed with strong homotopy actions
by $A$. A \emph{morphism of strong homotopy modules}
$f: L\ra M$ is a sequence of graded ($B$-linear) morphisms
\[
f_{i} : A\tp{i-1}\ten L \ra M
\]
homogeneous of degree $1-i$ such that for each $n\geq 1$, we
have
\begin{equation}\label{ainfmorph}
\begin{split}     
         & \sum _{k=1}^{n} 
    m_{k}(a_{1}, \dots , a_{k-1}, f_{n-k+1}(a_{k}, \dots , a_{n-1},x)) \\
\quad  = \quad & \;\; \sum _{i=1}^{n-2} (-1)^{i-1} 
         f_{n-1}(a_{1},\dots ,a_{i}a_{i+1},\dots ,a_{n-1},x) \\
    & +\sum _{k=1}^{n} (-1)^{n-k} 
f_{n-k+1}(a_{1},\dots ,a_{n-k}, m_{k}(a_{n-k+1},\dots ,a_{n-1},x)) \ko
\end{split}
\end{equation}
for all $a_i\in A$, $x\in L$.
For $n=1$, this specializes to
\[
m_{1}f_{1}=f_{1}m_{1}
\]
so that $f_{1}$ is a morphism of complexes. For $n=2$, we obtain
\begin{equation} \label{f1f2}
m_{1}(f_{2}(a_{1},x))+ m_{2}(a_{1},f_{1}(x)) =
f_{1}(m_{2}(a_{1},x))- f_{2}(a_{1},m_{1}(x)) \ko a_1\in A\ko x\in L\ko
\end{equation}
which means that for each $a_{1}\in A$, the morphism $f_{1}$ commutes 
with the left multiplication by $a_{1}$ up to the homotopy 
$x\mapsto f_{2}(a_{1},x)$. 

A morphism $f:L\ra M$ of strong homotopy modules is 
\emph{nullhomotopic} if there exists an \emph{homotopy from $f$
to $0$}, i.e. a family 
\[
h_{i}: A\tp{i-1}\ten L \ra  M\ko i\geq 1\ko 
\]
of graded morphisms homogeneous of degree $-i$ such that
for each $n\geq 1$, we have
\begin{equation}\label{ainfhomotopy}
\begin{split}     
f_{n} = & \sum _{k=1}^{n} 
    (-1)^{k-1} m_{k}(a_{1}, \dots , a_{k-1}, h_{n-k+1}(a_{k}, 
                                        \dots , a_{n-1},x)) \\
    & +\sum _{i=1}^{n-2} (-1)^{i-1} 
         h_{n-1}(a_{1},\dots ,a_{i}a_{i+1},\dots ,a_{n-1},x) \\
    & +\sum _{k=1}^{n} (-1)^{n-k} 
      h_{n-k+1}(a_{1},\dots ,a_{n-k}, m_{k}(a_{n-k+1},
                                          \dots ,a_{n-1},x)) \ko
\end{split}
\end{equation}
For $n=1$, this equation becomes 
\[
f_{1}= m_{1}\, h_{1}+ h_{1}\,m_{1}\ko 
\]
which means that $f_{1}$ is nullhomotopic. Two morphisms
between strong homotopy modules are \emph{homotopic} if their
difference is nullhomotopic.

We extend our dictionary
between strong homotopy and differential graded
coalgebra: Let $X=C\ten L$ and $Y=C\ten M$ be the differential 
graded comodules associated with $L$ and $M$, in analogy with
subsection \ref{differentialcoalgebra}. It is easy to check that
the map $f \mapsto \eps \circ f$, where $\eps :Y \ra M$ is the
canonical projection, is a bijection from the set of comodule morphisms
to the set of morphisms of graded modules $X \ra M$ and that under
this bijection, the morphisms of differential comodules correspond
exactly to the morphisms of strong homotopy modules. If
$f$ is a nullhomotopic morphism of differential graded
comodules, the map $h \mapsto \eps \circ h$ also induces
a bijection from the set of homotopies from $f$ to $0$ to
the set of homotopies from $\eps \circ f$ to $0$.

\section{From strong homotopy actions to strict actions}\label{strongstrict}

Let $k$ be a commutative ring and $A$ and $B$ (associative, unital) 
$k$-algebras. Let $L$ and $M$ be strong homotopy modules
(cf. \ref{morphismsstronghtpy}). We define a complex of $k$-modules
$\Hominfdot (L,M)$ as follows: its $p$th component is the $k$-module
of sequences of graded ($B$-linear) morphisms
\[
f_{n}: A\tp{n-1}\ten L\ra L\ko n\geq 1\ko 
\]
of degree $p+1-n$; its differential maps a sequence $f_{n}$
in $\Hominf ^{p}(L,M)$
to the sequence $g_{n}$ defined by  
\begin{equation}\label{ainfdifferential}
\begin{split}     
g_{n} = & \sum _{k=1}^{n} 
    (-1)^{p(k-1)} m_{k}(a_{1}, \dots , a_{k-1}, f_{n-k+1}(a_{k}, 
                                        \dots , a_{n-1},x)) \\
    & -(-1)^{p} \sum _{i=1}^{n-2} (-1)^{i-1} 
         f_{n-1}(a_{1},\dots ,a_{i}a_{i+1},\dots ,a_{n-1},x) \\
    & -(-1)^{p} \sum _{k=1}^{n} (-1)^{n-k} 
      f_{n-k+1}(a_{1},\dots ,a_{n-k}, m_{k}(a_{n-k+1},
                                          \dots ,a_{n-1},x)) \ko
\end{split}
\end{equation}
\begin{lemma}\label{ainfcomplex}
The square of the above differential vanishes. The group of
zero cycles of $\Hominfdot (L,M)$ identifies with the group of
morphisms of strong homotopy modules $L\ra M$. 
Zero boundaries correspond exactly to the nullhomotopic morphisms.
\end{lemma}

\begin{proof} We use the correspondence of \ref{morphismsstronghtpy}:
The set $\Hominf ^{p}(L,M)$ is in bijection with the
set of graded comodule morphisms $f:C\ten L \ra C\ten M$ of degree $p$
and the differential corresponds to the map
\[
f \mapsto b\circ f -(-1)^{p} f\circ b.
\]
This shows that we have a well-defined complex. The rest
follows upon inspection of \eqref{ainfmorph} and \eqref{ainfhomotopy}.
\end{proof}

Let $\Shmod{A}{B}$ denote the category of strong homotopy $A$-modules
over $B$ and let $\Bimod{A}{B}$ denote the category of complexes
of $A$-$B$-bimodules. We have an obvious functor 
\[
R: \Bimod{A}{B} \ra \Shmod{A}{B}
\]
which maps a complex of $A$-$B$-bimodules to the underlying 
$\Z$-graded $B$-module endowed with the homotopy action
given by the differential, the multiplication, and 
$m_{n}=0$ for all $n\geq 3$. We will construct a left adjoint.
We use the notations of \ref{differentialcoalgebra}.
For $X\in \Shmod{A}{B}$, let $LX$ be the complex whose underlying
graded $A$-$B$-bimodule is $A \ten C\ten X$ and whose differential
is
\begin{equation}
\begin{split}
d(a_{0},a_{1}, \dots , a_{n-1},x) = 
 & - (a_{0}a_{1}, a_{2}, \dots , a_{n-1},x) \\
 & + \sum _{i=1}^{n-2} 
     (-1)^{i-1}(a_{0}, \dots ,a_{i}a_{i+1}, \dots , a_{n-1},x) \\
 & + \sum _{k=1}^{n} (-1)^{n-k} 
     (a_{0}, a_{1}, \dots , a_{n-k}, m_{k}(a_{n-k+1}, \dots , a_{n-1},x)).
\end{split}
\end{equation}

\begin{lemma}\label{bimoduledifferential} The square of the above
differential vanishes.
\end{lemma}

\begin{proof}
Define a differential on the graded module $A\ten C$ by
\[
d(a_{0}, a_{1}, \dots ,a_{n-1})= - (a_{0}a_{1}, a_{2}, \dots, a_{n-1})
+ a_{0}\ten d_{C}(a_{1}, a_{2}, \dots , a_{n-1}).
\]
It is not hard to check that its square vanishes. On the other
hand, $C\ten X$ is endowed with the differential of \eqref{diffcomodule}.
Now the morphism
\[
\id _{A} \ten \Delta \ten \id _{X} : A\ten C \ten  X \ra (A\ten C)\ten (C\ten X)
\]
defines an isomorphism onto a graded submodule and the differential
given on $A\ten C\ten X$ is induced by the one on 
$(A \ten C)\ten (C\ten X)$.

\end{proof}

We define
\[
\phi : \Hominfdot (X,RY) \iso \Homdot (LX, Y)
\]
as the composition
\[
\Hominfdot (X,RY)=\Homdot_{k}(C\ten X,Y) \iso  
\Homdot_{A}(A\ten C\ten X,Y) =
\Homdot (LX,Y).
\]
More explicitly, a morphism of strong homotopy modules $f$ corresponds
to a morphism of graded modules $f: C\ten X \ra Y$. By definition,
the image of $f$ under $\phi $ maps $a\ten c \ten x$ to 
$a f(c,x)$. Of course, $\phi $ is an isomorphism of
graded $k$-modules. Its inverse
maps $g$ to $c\ten x \mapsto g(1\ten c\ten x)$. Note that
the fact that $A$ has a unit is crucial for this.

\begin{lemma}\label{adjointness} The isomorphism $\phi$ is
compatible with the differentials. In particular, the
functors $L$ and $R$ are adjoints and induce a pair
of adjoint functors in the homotopy categories.
\end{lemma}

\begin{proof}
The claim follows from the equalities
\begin{align*}
d(\phi (f)(a_{0}, \dots  & , a_{n-1},x)) \\
       & = a_{0} d(f(a_{1}, \dots ,a_{n-1},x)) \\
  \phi (f)(d(a_{0}, \dots & , a_{n-1},x)) \\
       &= -(a_{0} a_{1}) f(a_{2}, \dots , a_{n-1},x) \\
       &\quad + \sum _{i=1}^{n-2} (-1)^{i-1} 
          a_{0} f(a_{1}, \dots ,a_{i}a_{i+1}, \dots a_{n-1},x) \\
       &\quad + \sum _{k=1}^{n} (-1)^{n-k} 
          a_{0} f(a_{1},\dots ,a_{n-k},m_{k}(a_{n-k+1}, \dots , a_{n-1},x)) \\
\phi (d(f))(a_{0},\dots & ,a_{n-1},x) \\
       &= a_{0} d(f(a_{1},\dots ,a_{n-1},x))  
          -(-1)^{p} a_{0}(a_{1} f(a_{2},\dots a_{n-1},x)) \\
       &\quad - (-1)^p\sum _{i=1}^{n-2} (-1)^{i-1} 
          a_{0}f(a_{1},\dots, a_{i}a_{i+1},\dots a_{n-1},x) \\
       &\quad - (-1)^p\sum _{k=1}^{n} (-1)^{n-k} 
                            a_{0} f(a_{1}, \dots , a_{n-k},
          m_{k}(a_{n-k+1},\dots ,a_{n-1},x)).
\end{align*}
\end{proof}

The functors $R$ and $L$ induce a pair of adjoint functors
between the homotopy categories of $\Shmod{A}{B}$ and $\Bimod{A}{B}$.
A \emph{quasi-isomorphism} of $\Shmod{A}{B}$ is a morphism $f:X\ra Y$
such that $f_{1}$ is a quasi-isomorphism of the underlying complexes.
Then clearly the functor $R$ preserves quasi-isomorphisms. 
Since $A$ is projective over $k$, the functor $L$ also preserves
quasi-isomorphisms. Hence if we define the derived categories
$\cd \Shmod{A}{B}$ and $\cd \Bimod{A}{B}$ to be the localizations
of the homotopy categories
with respect to the quasi-isomorphisms, then $L$ and $R$ induce
a pair of adjoint functors between the derived categories:
\[
\begin{array}{c} \cd \Bimod{A}{B} \\
L \ua \da R \\
\cd \Shmod{A}{B}
\end{array}
\]
Let $Y$ be a complex of $A$-$B$-bimodules. It is easy to
check that $LRY = A\ten T(A[1])\ten Y$ is isomorphic to
the Hochschild resolution of $Y$ and that the 
adjunction morphism
\[
LRY = A\ten T(A[1]) \ten Y \ra Y
\]
identifies with the augmentation of the Hochschild resolution. 
In particular,
the adjunction is a quasi-isomorphism. It follows that the
functor $R: \cd \Bimod{A}{B}\ra \cd \Shmod{A}{B}$  is fully faithful.
If $X$ is a strong homotopy module, the adjunction morphism
$X \ra RL X = A \ten T(A[1]) \ten X$ is the morphism of strong
homotopy modules whose component in degree $i$ is the 
morphism 
\[
f_i: A\tp{i-1}\ten X \ra A \ten A[1]\tp{i-1}\ten X 
                             \subset A \ten T(A[1]) \ten X 
\]
given by
\[
(a_1, \ldots, a_{i-1})\ten x \mapsto 1 \ten (a_1, \ldots, a_{i-1})\ten x.
\]
We say that $X$ is \emph{$H$-unital} if the adjunction morphism
is a quasi-isomorphism. We deduce the

\begin{proposition}\label{equivainfstrict}
The functor
\[
R:\cd \Bimod{A}{B}\ra \cd \Shmod{A}{B}
\]
is an equivalence 
onto the full subcategory of $H$-unital strong homotopy modules.
Its inverse is induced by the functor $L$.
\end{proposition}

For the applications, we need a criterion for $H$-unitality:

\begin{lemma}\label{hunitality} Let $X$ be a strong homotopy module.
Then the following are equivalent:
\begin{itemize}
\item [(i)] $X$ is $H$-unital. 
\item [(ii)] The morphism of complexes of $k$-modules
      $m_{2}(1,?): X \ra X$ induces the identity in homology.
\end{itemize}
\end{lemma}

\begin{proof} 
Suppose that (i) holds. The square
\[
\begin{array}{rcl}
X & \ra & RL X \\
m_{2}(1,?)\da\; &   & \;\;\da m_{2}(1,?) \\
X & \ra & RL X
\end{array}
\]
commutes in the homotopy category of complexes of $k$-modules
thanks to \eqref{f1f2}.
By assumption, the adjunction morphism $X\ra RL X$ is a
quasi-isomorphism. The right vertical arrow is the identity
(since the $A$-module structure on $RLX= A \ten C \ten X$ 
is induced from that of $A$). So if we apply the homology 
functor to the diagram, we see that (ii) holds.

Suppose that (ii) holds. Consider the filtrations 
\[
F_{p}RL X = A\ten X \oplus (A\ten A[1]\ten X) \oplus \dots 
\oplus (A\ten A[1]\tp{p} \ten X) \ko p\geq 0\ko 
\]
and 
\[
F_{p}X=X\ko p\geq 0.
\]
The morphism
$f:X\ra RL X$, $x \mapsto 1\ten x$, is compatible with the filtrations.
The $E_{1}$-term of the spectral sequence associated
with $F_{p}RLX$ is the Hochschild resolution of the
graded $A$-module $\Hs X$. By our assumption, this module
is \emph{unital} and thus the Hochschild resolution is 
quasi-isomorphic to the module $\Hs X$ and the 
map $x \mapsto 1\ten x$ induces a quasi-isomorphism.
It follows that $f$ induces an isomorphism in the
$E_{2}$-terms of the spectral sequences. Since the
filtrations are bounded below and exhaustive,
the spectral sequences converge (by the classical convergence
theorem \cite[5.5.1]{Weibel94}) and $f$ is a
quasi-isomorphism.
\end{proof}

\section{Proof of existence}\label{proofexistence}

We prove part a) of theorem \ref{maintheorem}.
We put $m_1=d: T \ra T$ and construct $m_2$, $m_3$
as in section \ref{firstspecial}.
Since $A$ is projective over $k$, the vanishing
condition \eqref{Todastrongaction} follows from \eqref{Toda}.
Hence by theorem \ref{exstrongaction}, the triple $m_{1},m_{2},m_{3}$
may be completed to a strong homotopy action $m_{n},n\geq 1$, 
of $A$ on $T$. Let us denote by $\tilde{T}\in \Shmod{A}{B}$ the
corresponding strong homotopy module, cf. \ref{morphismsstronghtpy}.
In the homotopy action of $A$ on $T$, the unit of $A$ acts
by the identity, so that $m_{2}(1,?): \tilde{T}\ra \tilde{T}$ is
homotopic to the identity. By lemma \ref{hunitality}, the strong homotopy
module $\tilde{T}$ is $H$-unital. So by proposition \ref{equivainfstrict},
it comes from a complex of bimodules.
More precisely, the canonical morphism of strong homotopy modules
\[
f:\tilde{T} \ra RL \tilde{T}
\]
is a quasi-isomorphism. This means that 
$f_{1}$ is a quasi-isomorphism, which, by \eqref{f1f2}, is 
compatible with the homotopy actions of $A$ on $\tilde{T}$ and $L\tilde{T}$.
We put $X=L\tilde{T}= A \ten T(A[1]) \ten T$ and $\phi=f_1$. 

Note that if $T$ is a bounded complex, we 
can truncate $X$ to a bounded complex, as we did in
sections \ref{firstspecial} and \ref{secondspecial}.

\providecommand{\bysame}{\leavevmode\hbox to3em{\hrulefill}\thinspace}

\end{document}